\input amstex 
\documentstyle{amsppt}
\input bull-ppt
\keyedby{bull429/car}

\topmatter
\cvol{29}
\cvolyear{1993}
\cmonth{October}
\cyear{1993}
\cvolno{2}
\cpgs{235-242}
\title Stokes' theorem for nonsmooth chains\endtitle
\author Jenny Harrison\endauthor
\address  Department of Mathematics, University of 
California, Berkeley,
California 94720\endaddress
\ml harrison\@math.berkeley.edu\endml
\date February 1, 1993 and, in revised form, July 6, 
1993\enddate
\abstract Much of the vast literature on the integral 
during the last two
centuries concerns extending the class of integrable 
functions. In contrast,
our viewpoint is akin to that taken by Hassler Whitney 
[{\it Geometric
integration theory\/}, Princeton Univ. Press, Princeton, 
NJ, 1957] and by
geometric measure theorists because we extend the class of 
integrable {\it
domains\/}. Let $\omega$ be an $n$-form defined 
on $\Bbb R^m$. We show that if
$\omega$ is sufficiently smooth, it may be integrated over 
sufficiently
controlled, but nonsmooth, domains $\gamma$. The smoother 
is $
\omega$, the rougher may be $\gamma$. Allowable domains 
include a large class
of nonsmooth chains and topological $n$-manifolds immersed 
in $\Bbb R^m$. We
show that our integral extends the Lebesgue integral and 
satisfies a
generalized Stokes' theorem.\endabstract
\endtopmatter

\document

\heading 1. Introduction\endheading

The standard version of Stokes' theorem:
$$\int_{\partial M}\omega=\int _M\,d\omega$$
requires both a smooth $n$-manifold $M$ and a smooth 
$(n-1)$-form $\omega$. It
was realized at some point that one side of Stokes' 
formula could be used to
{\it define\/} the other side in more general situations. 
In particular Whitney 
\cite{14} used the right side to define the left in some 
cases where the
interior of $M$ is smooth, even though the boundary 
$\partial M$ is not smooth.
More generally the interior could be piecewise smooth or a 
suitable limit of
simplicial chains. Whitney systematically developed this 
insight by defining
his {\it flat norm\/} on chains, based on rectilinear 
subdivisions of simplices
of one higher dimension of which the chain in question is 
a partial boundary.
He treated forms as cochains---linear functions on the 
vector space of chains.

Stokes' theorem was thus extended by Whitney to 
integration of smooth forms
over objects that were limits in the flat norm of chains. 
These include certain
kinds of fractals. But other fractals escape Whitney's 
construction. Whitney's
example of a function nonconstant on a connected arc of 
critical points
\cite{15} shows the limits to any generalization of 
Stokes' theorem. Stokes' 
theorem (which for arcs is just the fundamental theorem of 
calculus) must 
fail for such arcs. There exists a Jordan curve in 3-space 
that has Hausdorff 
dimension $>2$, is not contained in any surface of finite 
area, and is not a
limit in the flat norm of simplicial  1-chains. Whitney's 
methods to not define
integration of forms over such a curve, nor do Lebesgue's. 
The methods for this
paper accomplish exactly this.

In a recursive construction we define a family of norms: 
For each real
$\lambda$ with $n\leq \lambda \leq m$ the {\it 
$\lambda$-natural norm\/}
$|A|^\natural_\lambda$ is defined on simplicial $n$-chains 
$A$ in $m$-space.
The definition utilizes projections of chains on 
hyperplanes, Whitney's method
of subdivisions and partial boundaries, and $\lambda$\<th 
powers of the
$n$-masses of the projected chains. The recursion is on 
the integer part of $
\lambda-n$.

The following theorem implies that this is indeed  a norm 
on oriented
simplicial chains in $\Bbb R^m$ (which is not obvious). 
Completion yields 
Banach spaces $X_{n,\lambda}=X_{n,\lambda}(\Bbb R^m)$. 
These spaces are nested,
containing increasingly complex chains as $\lambda$ 
increases. The parameter $
\lambda$ acts as a dimension; it bounds the fractal 
complexity of elements of
$X_{n,\lambda}$. The boundary operator $\partial$ on 
simplices extends to a 
continuous operator $\partial\:X_{n,\lambda}\rightarrow 
X_{n-1,\lambda}$.

\proclaim{Theorem 1} Integration of an $n$-form $\omega$ 
defines a bounded 
linear functional $L_\omega$ on $X_{n,\lambda}$ if 
$\omega$ is of class $C^{
\lambda-n}$. Furthermore, there is a constant $c>0$ such 
that for a simplicial
$n$-chain $A$\RM:
$$\left|\int_A\omega\right|\leq 
c|A|^\natural_\lambda\|\omega\|_{\lambda-n}$$
where $\|\omega\|_{\lambda-n}$ denotes the $\lambda-n$ 
norm on forms.
\endproclaim

This theorem extends results in \cite{14} and \cite7.

Denote the dual space of cochains by $X^{n,\lambda}$ and 
the product of a
cochain $X\in X^{n,\lambda}$ with a chain $A\in 
X_{n,\lambda}$ by $X\boldcdot
A$. Define $\Omega^{n,\lambda}=\{n$-forms $\omega$ of 
class $C^{\lambda-n}$
defined on $\Bbb R^m\}$. According to Theorem 1 the linear 
operator $A\mapsto 
\int_A\omega$ on simplicial $n$-chains in $\Bbb R^m$ is 
bounded and determines
an element $L_\omega\in X^{n,\lambda}$. This results in a 
linear mapping $L\:
\Omega^{n,\lambda}\rightarrow X^{n,\lambda}$. For $A\in 
X_{n,\lambda}$ we
define $\int_A\omega=L_\omega \boldcdot A$. Special cases: 
If $A$ is a
simplicial $n$-chain (finite) in $\Bbb R^m$, then $\int 
_A\omega$ is identical
to the Lebesgue integral for simplicial $n$-chains. If $A$ 
is an $n$-chain
(possibly infinite) in $\Bbb R^n$, then $\int_A\omega$ is 
identical to the
Lebesgue integral for infinite $n$-chains with finite 
$n$-mass. Since $L_{d
\omega}=dL_\omega$, we obtain the following: 

\proclaim{Theorem 2 \rm(Stokes' theorem for nonsmooth 
chains)} If $A\in X_{n,
\lambda}$ and $\omega\in C^{\lambda-n+1}$ is an 
$(n-1)$-form defined on $\Bbb
R^m$, then $\int_A\,d\omega=\int_{\partial 
A}\omega$.\endproclaim

\subheading{Historical comments}
In 1982 the author found a $C^{2+\alpha}$ counterexample 
to the Seifert 
Conjecture \cite {2, 3}. At the heart of her construction 
is a diffeomorphism
of the two-sphere with a fractal equator $\gamma$ as an 
invariant set.  The 
Hausdorff dimension of $\gamma$ is precisely $1+\alpha$.  
Coarser relations
between the differentiability class of a function and the 
topological dimension 
of a related set were well known in the theories of Sard 
\cite {13}, Denjoy
\cite {1}, and Whitney \cite {15}. Subsequently, these too 
were found to have
fractal versions.  (See \cite {6, 11}.) In each of these 
theories appears the
same phenomenon---the smoother the function, the rougher 
an associated set. 
This paper is a result of the author's search for a 
general principle
underlying this duality between dimension and 
differentiability.

Fractals are rife in many fields: geometric measure 
theory, dynamical systems,
PDEs, function theory, foliations\<$\dotsc$\<. Except for 
various dimension
theories, however, there are few available techniques for 
investigating them. 
Usually there is a large measure of geometric control over 
their  formation
that can be proved or reasonably postulated.  In at least 
some interesting
cases this ought to allow integration of sufficiently 
smooth forms over the
fractals.  For example, it would be very useful, and 
perhaps not too difficult,
to show that certain solenoid-like invariant sets for 
flows are domains for
integration of 1-forms and that certain closed unions of 
leaves of a foliation
by $k$-manifolds are domains for integration of
$k$-forms.  Such sets could then be proved fairly directly 
to carry homology
cases.  Currently such proofs (due in various contexts to 
S. Schwartzmann, J.
Plante, D. Sullivan, and others) are indirect and hard to 
come by.  It may well
be that these new discoveries in geometric integration 
theory will provide
simple but powerful new tools to sharpen differentiability
results and reduce technicalities in known proofs.

\heading 2. Norms on chains\endheading

We take Whitney's approach in \cite{14} with simplicial 
$n$-chains. It provides
an algebraic way of equating chains with common simplicial 
subdivisions. An
{\it $n$-simplex\/} $\sigma=p_0\cdots p_n$ in $\Bbb R^n$ 
is the convex hull of
$n+1$ vertices $p_i$ in $\Bbb R^n$, $i=0,\dots, n$. If the 
$n$ vectors
$\{p_0-p_i\}$ are linearly dependent, then $\sigma$ is 
{\it degenerate\/}. The
order of the vertices $p_i$ of a nondegenerate simplex 
$\sigma$ determines an
orientation on $\sigma$.  The simplex $-\sigma$ is 
identified with the  simplex
with the same
pointset as $\sigma$ but the opposite orientation. Define 
a {\it simplicial
$n$-chain\/} $A$ {\it in\/} $\Bbb R^n$ as an equivalence 
class of formal sums $
\sum a_i\sigma_i $ with real coefficients and oriented 
simplices in $\Bbb
R^n$ as follows: A formal  sum $ \sum a_i\sigma_i$ defines 
a function $A\:
\Bbb R^n\rightarrow \Bbb R$ by $A(x)= \sum \pm a_i$ for 
$x\in \operatorname{
int}(\sigma_i)$
where $+a_i$ is used if $\sigma_i$ is oriented in the same 
way as $\Bbb R^n$
and $-a_i$ is used otherwise. Set $A(x)=0$ if $x$ does not 
lie in the interior
of any simplex $\sigma_i$. Two formal sums of oriented 
simplices are equivalent
if the functions defined by them are identical except in 
sets with 
$n$-dimensional 
Lebesgue measure zero. (Degenerate simplices can be used 
but are
equated with the zero $n$-chain.) The definition of a 
simplicial $n$-chain $A$
in $\Bbb R^n$ is independent of the orientation chosen for 
$\Bbb R^n$.

More generally we can define a {\it simplicial $n$-chain 
in $\Bbb R^m$\/} when
$m\geq n$. For each affine $n$-plane $P$ in $\Bbb R^m$ let 
$C_n(P)$ denote the 
linear space of simplicial $n$-chains in $P$. A {\it 
simplicial $n$-chain $A$
in $\Bbb R^m$\/} is an element of the direct sum of the 
$C_n(P)$. Let $A$ and
$B$ be two simplicial $n$-chains in $\Bbb R^m$. If the 
corresponding summands
of $A$ and $B$ are equivalent, we write $A\sim B$.

Let $M_n(\sigma)$ denote the $n$-dimensional {\it mass\/} 
or volume of an
$n$-simplex $\sigma$. Let $A$ be a simplicial $n$-chain in 
$\Bbb R^m$ and $\pi$
the orthogonal projection onto an affine subspace of $\Bbb 
R^m$. For $0\leq 
\lambda \leq n$ define the {\it projected 
$\lambda$-mass\/} of $A$ to be
$$M_{\lambda,\pi}(A)=\inf\left\{\sum|a_i|(M_n(\pi 
\sigma_i))^{\lambda/n}:\sum
a_i\sigma_i\sim A\right\}\.$$
If $\pi$ is the identity, we write 
$M_\lambda(A)=M_{\lambda,\pi }(A)$ for 
simplicity of notation. Then $M_n(A)$ is the $n$-mass of 
$A$.
It is a norm on simplicial $n$-chains in $\Bbb R^m$.

\subheading{Integration of smooth integrands over smooth 
domains} The mass norm
$M_n$ is often used to estimate the integral of a 
continuous form over a domain with
finite mass. For example, let $A$ be a simplicial 
$n$-chain in $\Bbb R^m$ and $
\omega$ a continuous $n$-form. Clearly,
$$\left| \int_A \omega\right|\leq M_n(A)\|\omega\|_0\.$$
Whitney used a smaller norm on chains to find better 
estimates for 
simplicial $n$-chains with large mass.
\subheading{Whitney's flat norm} Let $A$ be a simplicial 
$n$-chain in $\Bbb
R^m$. Whitney defined the {\it flat\/} norm of $A$ as
$$|A|^\flat=\inf\{M_n(A-\partial C )+M_{n+1}(C)\},$$
where the infimum is taken over all simplicial $(n+
1)$-chains $C$. For example,
if $A$ is a Jordan curve in $\Bbb R^2$, $|A|^\flat$ is 
bounded by
the area of the 
Jordan domain spanning $A$. If $A$ is an arc in $\Bbb R^2$,
consider any simplicial 1-chain $B$ so that $A+B$ is a 
Jordan curve. $
|A|^\flat$ is bounded by the sum of the arc length of $B$ 
plus the area of the
Jordan domain spanning $A+B$.

The flat norm is used to complete the space of simplicial 
$n$-chains in $
\Bbb R^m$. Denote the resulting Banach space of {\it flat 
$n$-chains\/} by
$X_n$ and the dual space of {\it flat $n$-cochains\/} by 
$X^n$. Whitney showed
that the boundary operator $\partial$ defined on 
simplicial $n$-chains is 
continuous with respect to the flat norm and extends to 
$\partial\:X_n
\rightarrow X_{n-1}$.

Let $A$ be a simplicial $n$-chain and $C$ a simplicial $(n+
1)$-chain in $
\Bbb R^m$. Let $\omega$ be an $n$-form of class $C^1$ 
defined on $\Bbb R^m$.
Whitney used Stokes' theorem for simplicial $n$-chains to 
show 
$$\left|\int_A\omega \right|\leq \left|\int_{A-\partial 
C}\omega\right|+
\left|\int_C\,d\omega\right|\leq M_n(A-\partial 
C)\|\omega\|_0+M_{n+1}(C)\|d
\omega\|_0\.\tag"{$(*)$}"$$
Therefore, $|\int_A\omega|\leq |A|^\flat\|\omega\|_1$. 
This estimate enabled
Whitney to integrate over elements of $X_n$. Whitney 
believed he could
generalize Stokes' theorem and wrote at the end of his 
first paper of
geometric integration theory \cite{16}, ``With the help of 
methods described
above a very general form of Stokes' theorem may be 
proved. We shall not give
details here.'' However, it turned out that his integral 
does not extend to the
Lebesgue integral. Theorem 1 generalizes Whitney's theorem 
and leads to
strict generalizations of the Lebesgue integral and 
Stokes' theorem.

\subheading{The natural norms \cite5\rm} Let $\gamma$ be a 
Jordan curve in $\Bbb
R^3$ that has Hausdorff dimension greater than two. This 
curve not only has
infinite arc length but has no spanning surface with 
finite area. Thus a
sufficiently close simplicial 1-chain $A$ will only bound 
2-chains with huge
area. Neither the arc length norm $M_1(A)$ nor the flat 
norm $|A|^\flat$ will give
good estimates for integrals over $A$. The natural norms 
defined below give
sharper estimates.

Since $\partial \partial=0$, only simplicial $n$-chains 
without boundary have
spanning $(n+1)$-chains. Yet all simplicial $(n+1)$-chains 
$C$ can be viewed as
``partial'' spanning sets for an arbitrary simplicial 
$n$-chain $A$. The natural
norms take into account the $(n+1)$-mass of these partial 
spanning sets $C$ and
the projected $n$-mass of what is left over, namely, 
$A-\partial C$. (See 
Figure 1 on page 240.)

List the $n$-dimensional coordinate planes of $\Bbb R^m$ 
from 1 to $N=\binom
mn$. Let $\pi_i$ be the projection onto the $i$\<th 
coordinate plane in this
list.

Let $A$ be an $n$-chain in $\Bbb R^m$. If $0\leq \lambda 
\leq n$, define the
{\it $\lambda$-natural norm\/} of $A$ by
$$|A|^\natural_\lambda= \sum^N_{i=1}M_{\lambda ,\pi_i}(A 
)\.$$
Then $|A|^\natural_n$ is proportional to the $n$-mass 
$M_n(A)$.

For $n<\lambda\leq m$, $\lambda\in\Bbb R$,
 define the $\lambda$-natural norm of $A$ by recursion on
the integer part of $\lambda-n$:
$$|A|^\natural_\lambda= 
\sum^N_{i=1}\inf_C\{M_{n,\pi_i}(A-\partial
C)+|C|^\natural_\lambda\}\.$$
Each term is infimized over all simplicial $(n+1)$-chains 
$C$ in $\Bbb R^m$.
The flat norm corresponds to $\lambda=n+1$.  It is 
interesting to observe
that $\lambda$ may be an integer much larger than $n+1$. 
This allows us to work
with domains that are curves in $\Bbb R^3$, for example.  
The $\lambda$-natural
norm depends on the parameter $\lambda$, the topological 
dimension $n$ of $A$,
and the minimal ambient space dimension $m$ of $A$.  
However, if $A$ is an
$n$-chain in $\Bbb R^p$ and in $\Bbb R^{p'}$, the norms 
defined are identical.  
For simplicity of notation the dependency on $n$ and $m$ 
are suppressed.  The 
definition is independent of choice of coordinates.  That 
is, under a
$C^{\lambda-n+1}$ change of coordinates with bounded 
derivative, the
$\lambda$-natural norm changes to an equivalent norm.

As does Whitney, we complete and obtain a Banach space 
$X_{n,\lambda}$ of {\it $
\lambda$-natural $n$-chains\/} and its dual space 
$X^{n,\lambda}$ of {\it $
\lambda$-natural $n$-cochains\/}. 

\rem{Remark} There is a one-parameter family of homology 
and cohomology groups 
associated with the $\lambda$-natural norm, 
$H_{n,\lambda}$ and
$H^{n,\lambda}$.\endrem

Now let $\omega$ be a 1-form defined on $\Bbb R^m$ of 
class $C^1$ and $A\in 
X_{1,\lambda}$. We have seen that the integral of $\omega$ 
over $A$ can be
estimated using either arc length or the flat norm, both 
of which may be very
large.  Theorem 1 shows that if $\omega$ is of class 
$C^2$, the integral
can also be estimated using the 3-natural norm which may 
be much smaller.
$$\left| \int_A \omega\right|\leq 
|A|^\natural_3\|\omega\|_z\tag"{$(**)$}"$$

To prove Theorem 1,  we iterate a
modified version of the Stokes' argument as seen in $(*)$ 
(see \cite5 for
details; also see the example below). 
The reader should take note that Stokes' theorem, as such, 
may be
applied only once since $dd\omega=0$ and $\partial 
\partial A=0$. However, the 
components of $d\omega$ may not be exact, and we may apply 
the exterior
derivative to them. This means we may continue to use the 
exterior derivative
and iterate the ``Stokes' process'' in a nontrivial 
fashion for as many times
as $\omega$ is differentiable. We also use partial 
spanning sets of spanning
sets. 

Typically, an integral $\int_A\omega$ is treated with two 
basic methods: If $A=
\partial B$ and $\omega$ is smooth, one may ``go forward'' 
and integrate $d
\omega$ over $B$. Alternatively, if the form 
$\omega=d\nu$, one may ``go
backward'' and integrate $\nu$ over $\partial A$. 
Sometimes neither of the
equivalent new integrals is easier to calculate. Our 
methods allow one to go
forward many times under suitable conditions, taking the 
exterior derivatives
of forms in the hopes of finding an equivalent integral 
that may be calculated.
It is worth mentioning that it is also possible, under 
suitable conditions, to
go backward, taking antiderivatives of forms many times, 
and find an equivalent
integral that may be calculated.

\fg{10pc}
\caption{Figure 1 }

\ex{Example} There exists a self-similar
Jordan curve $\gamma$ in the unit cube with Hausdorff 
dimension $>2$ so that each projection $\pi_i(\gamma)$, 
$i=1,2,3$, onto the
$i$\<th two-dimensional coordinate plane is an immersed 
curve bounding an
infinite 2-chain with finite 2-mass $<1$. 
(The curve $\gamma$ is constructed by adding homothetic 
replicas of the curve
$\alpha $ to itself in Figure 1 so that $\gamma$ is 
embedded and has Hausdorff 
dimension $\sim\log 11/\log 3>2.)$
A simplicial 1-chain $A$ can be found
approximating $\gamma$ along with simplicial 2-chains 
$C_i$ spanning $A$,
$i=1,2,3$, such that $M_{2,\pi_i}(C_i)<1$. 
(Each surface $C_i$ is a sum of homothetic replicas of the 
surfaces $\sigma$
and $\tau$ in Figure 1, chosen to minimize the projected 
2-mass of $C_i$.)
Furthermore, the simplicial 3-chains
$D_{ij}$ spanning $C_i+C_j$ have $M_3(D_{ij})<1$. Even 
though the norms
$M_1(A)$ and $|A|^\flat$ are both large, the 3-natural 
norm $|A|^\natural_3$ 
is not, for
$$|C^\natural_i|_3\leq \sum_{j=1}^3(M_{2,\pi_j}(C_j)+
|D_{ij}|_3^\natural
)<12\quad \text{
and thus}\quad |A|^\natural_3\leq 
\sum^3_{i=1}|C_i|^\natural_3<36\.$$
This curve $\gamma$ cannot be a domain for the classical 
Lebesgue integral and
cannot be treated using Whitney's methods. We show in \S3 
how to choose $A_k$ 
approximating $\gamma$ and define $\int_\gamma \omega=\lim 
\int_{A_k}\omega$.

We verify that Theorem 1 (see $(**))$ is valid for this 
example. It suffices
to estimate
$|\int_A\omega_i|=|\int_{C_i}\,d\omega_i|$ where 
$\omega_i$ is a component of $
\omega$. Note that $C_i+C_j=\partial D_{ij}$. Thus
$$\left|\int_{C_i}\,d\omega_i\right|\leq \sum_{
j=1}^3\left|\int_{C_j}(d\omega_i)_j\right|+
\left|\int_{D_{ij}}\,d(d\omega_i)_j\right|
\leq |A|^\natural_3\|\omega\|_2\.$$\endex

\subheading{Relation to the Lebesgue integral} Our 
integral of smooth forms
contains as a special case the Lebesgue integral: First 
associate an $L^1$ 
function $f\:\Bbb R^1\rightarrow \Bbb R^1$ with an element 
of $X_{1,2}$ as
follows.
Assume for simplicity that $f$ is the 
characteristic function of a bounded, measurable set $E$. 
Let $S$ denote the
region under the graph $\Gamma$ of $f$ and above the 
$x$-axis, 
oriented positively. Let $m_n$ denote
$n$-dimensional Lebesgue measure. Since $f
\in L^1$, $m_2(S)<\infty $. For each $k$ there exists a 
union $P_k$ of disjoint
intervals such that $E\subset  P_k$ and $m_1(P_k)<\infty 
$. Furthermore,
$\bigcap P_k=E$ except for a zero set. Subdivide each of 
the rectangles of $P_k
\times I$ into finitely many 2-simplices and orient each 
positively. Denote the
formal sum of these oriented 2-simplices by $Q_k$. Since 
$m_2(P_k\times
I)<\infty $, each $Q_k$ is an element of $X_{2,2}$. Since 
$m_2(S)<\infty $, the 
sequence $\{Q_k\}$ is Cauchy and thus converges to $A_f\in 
X_{2,2}$. It can be
shown that $\partial A_f\in X_{1,2}$ is canonically 
associated with $f$. 
Now it is well known (see Saks \cite{12}) that the 
Lebesgue integrals $\int_{
\Bbb R}f$ and $\int_S\,dx\,dy$ are equal. Since $S$ is 
canonically represented
by $A_f\in X_{2,2}$, the Lebesgue integral 
$\int_S\,dx\,dy$ equals our integral
$\int_{A_f}\,dx\,dy$. Since $ydx$ is analytic and 
$\partial A_f\in X_{1,2}$, 
Theorem 2 implies the Lebesgue integral $\int_{\Bbb R}f$ 
equals our integral
$\int_\Gamma y\,dx$. (This example can be treated with the 
flat norm alone.)

\heading 3. Domains of integration\endheading

Which oriented $n$-dimensional topological submanifolds $M 
\subset \Bbb R^m$ 
with  
boundary can serve as domains for integrals of smooth 
forms $\omega$
of class $C^{\lambda-n}$, $n\leq \lambda \leq m$?
The idea is to choose
appropriate simplicial approximators $A_k$ for $M$ that 
converge in $X_{n,
\lambda}$ for $\lambda$ sufficiently large and define 
$\int_M\omega=\lim
\int_{A_k}\omega$. (See \cite4 for some  numerical 
methods.) But choosing $A_k$
takes care. For example, let $M_1$ be a Jordan curve in 
$\Bbb R^2$ with
positive Lebesgue area. Any Cauchy sequence of simplicial 
1-chains inside
$M_1$ will have different limit point in $X_{1,2}$ from a 
Cauchy sequence of
1-chains outside $M_1$. The area of $M_1$ itself 
contributes an unavoidable
error. 

Consider a compact oriented $n$-manifold $M\subset  \Bbb 
R^m$ with boundary.
There is an algorithm for constructing a particular 
sequence of simplicial
$n$-chains in the nerves of coverings of $M$ by boxes in 
cubic lattices, the
inverse limit of these chains representing the fundamental 
\v Cech homology
class of $(M,\partial M)$. If this sequence of {\it binary 
approximators\/}
converges to some $\psi(M)\in X_{n,\lambda}(\Bbb R^m)$, 
then $M$ is a
$(\lambda,n)$-{\it set\/}. The algorithm commutes with the 
boundary operator on
chains. Examples of $(\lambda,n)$-sets include some planar 
Jordan curves with 
positive Lebesgue area and the above example of a Jordan 
curve in 3-space with
Hausdorff dimension $>2$.
If $M$ is a $(\lambda,n)$-set and $\omega$ is an
$n$-form of class $C^{\lambda-n}$, we can define 
$\int_M\omega=\int_{\psi(M)}
\omega$ and apply Theorem 2 to conclude:
\proclaim{Theorem 3 \rm (Stokes' theorem for 
$(\lambda,n)$-sets)} If $M$ is a $(
\lambda,n)$-set and $\omega\in C^{\lambda-n+1}$, then
$$\int_{\partial M}\omega=\int_M\,d\omega\.$$
\endproclaim

Theorem 3 extends the classic Green's theorem for 
codimension one boundaries $
\gamma$ and $C^1$ forms $\omega$. H\"older versions of the 
classic Green's 
theorem relating the box dimension of $\gamma$ to the 
H\"older exponent of $
\omega$ are proved in \cite{7, 8, 10}.

There are many examples of $(\lambda,n)$-sets:

(i) There exist Jordan curves in $\Bbb R^n$ with Hausdorff 
dimension $d$ that
are $(d,1)$-sets for every $1\leq d\leq n$.

(ii) Every hypersurface immersed in $
\Bbb R^m$ with Hausdorff dimension $<m$ is an $(m,m-1)$-set.

(iii) Compact oriented
$n$-manifolds immersed in $\Bbb R^m$ that are locally 
Lipschitz
graphs are $(n,n)$-sets.

(iv) Compact oriented 1-manifolds immersed in $\Bbb R^2$ 
that are locally graphs of
H\"older functions with exponent $\alpha $ are $(2/(1+
\alpha ),1)$-sets.

\rem{Remark} Compactness in examples (ii), (iii), and (iv) 
can be relaxed.
An arc in $\Bbb R^2$ that spirals to the origin passing 
through the $x$-axis at
$x_n=1/\sqrt n$, $n\geq1$, is {\it not\/} a 
$(\lambda,1)$-set for any $1\leq 
\lambda \leq 2$. Other domains of integration include 
boundaries of open sets
such as the topologist's since circle and the Denjoy 
Cantor set.\endrem

\subheading{Integration on manifolds}
We pass to smooth manifolds by way of local coordinates, 
since the
$\lambda$-natural norm changes to an equivalent norm under a
$(\lambda-n+1)$-smooth change of coordinates with bounded 
derivative. One must
also extend the definition of these norms to chains in 
proper open subsets of $
\Bbb R^m$ as Whitney did for his flat norm. The integral 
can then be defined
for smooth forms defined on smooth compact $m$-manifolds 
$M$ over domains that
are ``rough'' subsets of $M$.

\heading Acknowledgment \endheading 
The author thanks Morris Hirsch for his helpful comments 
and detailed
verification of the main results.

\enddocument